\documentclass[oneside,12pt]{article}
\usepackage{hyperref} 
\usepackage{tocloft}                   
\usepackage{indentfirst}
\usepackage{amsmath}
\usepackage{amsfonts}
\usepackage{amssymb}
\usepackage{pict2e}
\usepackage{stmaryrd}
\usepackage[x11names,svgnames]{xcolor} 
\usepackage{proof}   
\input diagxy        
\xyoption{curve}     
\usepackage{edrx15}               

\makeatletter
\ifluatex
  \usepackage{mynewunicodechar}   
\else
  \def\defunicodechar#1#2{\@namedef{u8:\detokenize{#1}}{#2}}
\fi
\makeatother

\ifluatex

\defunicodechar{á}{\'a}
\defunicodechar{é}{\'e}
\defunicodechar{í}{\'\i}
\defunicodechar{ó}{\'o}
\defunicodechar{ú}{\'u}
\defunicodechar{Á}{\'A}
\defunicodechar{É}{\'E}
\defunicodechar{Í}{\'I}
\defunicodechar{Ó}{\'O}
\defunicodechar{Ú}{\'U}

\defunicodechar{à}{\`a}
\defunicodechar{è}{\`e}
\defunicodechar{ì}{\`\i}
\defunicodechar{ò}{\`o}
\defunicodechar{ù}{\`u}
\defunicodechar{À}{\`A}
\defunicodechar{È}{\`E}
\defunicodechar{Ì}{\`\I}
\defunicodechar{Ò}{\`O}
\defunicodechar{Ù}{\`U}

\defunicodechar{â}{\^a}
\defunicodechar{ê}{\^e}
\defunicodechar{î}{\^\i}
\defunicodechar{ô}{\^o}
\defunicodechar{û}{\^u}
\defunicodechar{Â}{\^A}
\defunicodechar{Ê}{\^E}
\defunicodechar{Î}{\^\I}
\defunicodechar{Ô}{\^O}
\defunicodechar{Û}{\^U}

\defunicodechar{ã}{\~a}
\defunicodechar{õ}{\~o}
\defunicodechar{Ã}{\~A}
\defunicodechar{Õ}{\~O}

\defunicodechar{ç}{\c c}
\defunicodechar{Ç}{\c C}

\defunicodechar{°}{^\circ}
\defunicodechar{¹}{^{-1}}
\defunicodechar{²}{^2}
\defunicodechar{³}{^3}
\defunicodechar{±}{\pm}
\defunicodechar{÷}{\div}
\defunicodechar{·}{\cdot}
\defunicodechar{×}{\times}
\defunicodechar{¬}{\neg}
\defunicodechar{§}{\S}
\defunicodechar{¡}{\text{\textexclamdown}}

\defunicodechar{Δ}{\Delta}
\defunicodechar{Γ}{\Gamma}
\defunicodechar{Θ}{\Theta}
\defunicodechar{Π}{\Pi}
\defunicodechar{Λ}{\Lambda}
\defunicodechar{Σ}{\Sigma}
\defunicodechar{Φ}{\Phi}
\defunicodechar{Ψ}{\Psi}
\defunicodechar{Ω}{\Omega}

\defunicodechar{α}{\alpha}
\defunicodechar{β}{\beta}
\defunicodechar{γ}{\gamma}
\defunicodechar{δ}{\delta}
\defunicodechar{ε}{\epsilon}
\defunicodechar{ζ}{\zeta}
\defunicodechar{η}{\eta}
\defunicodechar{θ}{\theta}
\defunicodechar{ι}{\iota}
\defunicodechar{κ}{\kappa}
\defunicodechar{λ}{\lambda}
\defunicodechar{μ}{\mu}
\defunicodechar{ν}{\nu}
\defunicodechar{ξ}{\xi}
\defunicodechar{π}{\pi}
\defunicodechar{ρ}{\rho}
\defunicodechar{σ}{\sigma}
\defunicodechar{τ}{\tau}
\defunicodechar{φ}{\phi}
\defunicodechar{ϕ}{\origphi}
\defunicodechar{χ}{\chi}
\defunicodechar{ψ}{\psi}
\defunicodechar{ω}{\omega}

\defunicodechar{•}{\bullet}
\defunicodechar{⅋}{\bindnasrepma}
\defunicodechar{←}{\ot}
\defunicodechar{↑}{\upto}
\defunicodechar{→}{\to}
\defunicodechar{↓}{\dnto}
\defunicodechar{↔}{\bij}
\defunicodechar{↕}{\updownarrow}
\defunicodechar{↖}{\nwarrow}
\defunicodechar{↗}{\nearrow}
\defunicodechar{↘}{\searrow}
\defunicodechar{↙}{\swarrow}
\defunicodechar{↣}{\epito}
\defunicodechar{↣}{\twoheadrightarrow}
\defunicodechar{↤}{\mapsfrom}
\defunicodechar{↦}{\mapsto}
\defunicodechar{⇐}{\Leftarrow}
\defunicodechar{⇒}{\funto}
\defunicodechar{∀}{\forall}
\defunicodechar{∂}{\partial}
\defunicodechar{∃}{\exists}
\defunicodechar{∅}{\emptyset}
\defunicodechar{∇}{\nabla}
\defunicodechar{∈}{\in}
\defunicodechar{∖}{\backslash}
\defunicodechar{∘}{\circ}
\defunicodechar{√}{\sqrt}
\defunicodechar{∞}{\infty}
\defunicodechar{∧}{\land}
\defunicodechar{∨}{\lor}
\defunicodechar{∩}{\cap}
\defunicodechar{∩}{\cap}
\defunicodechar{∪}{\cup}
\defunicodechar{≃}{\simeq}
\defunicodechar{≅}{\cong}
\defunicodechar{≠}{\neq}
\defunicodechar{≤}{\le}
\defunicodechar{≥}{\ge}
\defunicodechar{⊂}{\subset}
\defunicodechar{⊃}{\supset}
\defunicodechar{⊆}{\subseteq}
\defunicodechar{⊇}{\supseteq}
\defunicodechar{⊓}{\sqcap}
\defunicodechar{⊔}{\sqcup}
\defunicodechar{⊢}{\vdash}
\defunicodechar{⊨}{\vDash}
\defunicodechar{⊣}{\dashv}
\defunicodechar{⊤}{\top}
\defunicodechar{⊥}{\bot}
\defunicodechar{⋀}{\bigwedge}
\defunicodechar{⋁}{\bigvee}
\defunicodechar{⋂}{\bigcap}
\defunicodechar{⋃}{\bigcup}
\defunicodechar{♭}{\flat}
\defunicodechar{♮}{\natural}
\defunicodechar{♯}{\sharp}
\defunicodechar{⟦}{\llbracket}
\defunicodechar{⟧}{\rrbracket}
\defunicodechar{⠆}{{:}}

\defunicodechar{∼}{\sim}
\defunicodechar{≈}{\approx}
\defunicodechar{≡}{\equiv}

\defunicodechar{⋄}{\lozenge}
\defunicodechar{⋅}{\Box}
\defunicodechar{◻}{\Box}
\defunicodechar{□}{\Box}
\defunicodechar{…}{\ldots}
\defunicodechar{∫}{\int}
\defunicodechar{⇀}{\rightharpoonup}
\defunicodechar{〈}{\langle}
\defunicodechar{〉}{\rangle}
\defunicodechar{⊙}{\odot}
\defunicodechar{⊖}{\ominus}
\defunicodechar{⊕}{\oplus}
\defunicodechar{⊘}{\oslash}
\defunicodechar{⊗}{\otimes}
\defunicodechar{▁}{\_}

\defunicodechar{𝐛}{\mathbf}
\defunicodechar{𝐢}{\textsl}
\defunicodechar{𝐫}{\mathrm}
\defunicodechar{𝐭}{\text}
\defunicodechar{𝐬}{\mathsf}

\else

\usepackage[utf8]{inputenc}

\DeclareUnicodeCharacter{00B0}{^\circ}             
\DeclareUnicodeCharacter{00B9}{^{-1}}              
\DeclareUnicodeCharacter{00B2}{^2}                 
\DeclareUnicodeCharacter{00B3}{^3}                 
\DeclareUnicodeCharacter{00B1}{\pm}                
\DeclareUnicodeCharacter{00F7}{\div}               
\DeclareUnicodeCharacter{00B7}{\cdot}              
\DeclareUnicodeCharacter{00D7}{\times}             
\DeclareUnicodeCharacter{00AC}{\neg}               
\DeclareUnicodeCharacter{00A7}{\S}                 
\DeclareUnicodeCharacter{00A1}{\text{\textexclamdown}}
\DeclareUnicodeCharacter{03BC}{\mu}                
\DeclareUnicodeCharacter{03BD}{\nu}                
\DeclareUnicodeCharacter{03C1}{\rho}               
\DeclareUnicodeCharacter{03C3}{\sigma}             
\DeclareUnicodeCharacter{03C4}{\tau}               
\DeclareUnicodeCharacter{03C6}{\phi}               
\DeclareUnicodeCharacter{03C7}{\chi}               
\DeclareUnicodeCharacter{03C8}{\psi}               
\DeclareUnicodeCharacter{0394}{\Delta}             
\DeclareUnicodeCharacter{0393}{\Gamma}             
\DeclareUnicodeCharacter{0398}{\Theta}             
\DeclareUnicodeCharacter{03A0}{\Pi}                
\DeclareUnicodeCharacter{03A3}{\Sigma}             
\DeclareUnicodeCharacter{03A6}{\Phi}               
\DeclareUnicodeCharacter{03A8}{\Psi}               
\DeclareUnicodeCharacter{03A9}{\Omega}             
\DeclareUnicodeCharacter{03B1}{\alpha}             
\DeclareUnicodeCharacter{03B2}{\beta}              
\DeclareUnicodeCharacter{03B3}{\gamma}             
\DeclareUnicodeCharacter{03B4}{\delta}             
\DeclareUnicodeCharacter{03B5}{\epsilon}           
\DeclareUnicodeCharacter{03B7}{\eta}               
\DeclareUnicodeCharacter{03B8}{\theta}             
\DeclareUnicodeCharacter{03B9}{\iota}              
\DeclareUnicodeCharacter{03BB}{\lambda}            
\DeclareUnicodeCharacter{03C0}{\pi}                
\DeclareUnicodeCharacter{03C9}{\omega}             
\DeclareUnicodeCharacter{03D5}{\origphi}           
\DeclareUnicodeCharacter{2022}{\bullet}            
\DeclareUnicodeCharacter{214B}{\bindnasrepma}      
\DeclareUnicodeCharacter{2190}{\ot}                
\DeclareUnicodeCharacter{2191}{\upto}              
\DeclareUnicodeCharacter{2192}{\to}                
\DeclareUnicodeCharacter{2193}{\dnto}              
\DeclareUnicodeCharacter{2194}{\bij}               
\DeclareUnicodeCharacter{2195}{\updownarrow}       
\DeclareUnicodeCharacter{2196}{\nwarrow}           
\DeclareUnicodeCharacter{2197}{\nearrow}           
\DeclareUnicodeCharacter{2198}{\searrow}           
\DeclareUnicodeCharacter{2199}{\swarrow}           
\DeclareUnicodeCharacter{21A3}{\epito}             
\DeclareUnicodeCharacter{21A3}{\twoheadrightarrow} 
\DeclareUnicodeCharacter{21A4}{\mapsfrom}          
\DeclareUnicodeCharacter{21A6}{\mapsto}            
\DeclareUnicodeCharacter{21D0}{\Leftarrow}         
\DeclareUnicodeCharacter{21D2}{\funto}             
\DeclareUnicodeCharacter{2200}{\forall}            
\DeclareUnicodeCharacter{2203}{\exists}            
\DeclareUnicodeCharacter{2205}{\emptyset}          
\DeclareUnicodeCharacter{2208}{\in}                
\DeclareUnicodeCharacter{2216}{\backslash}         
\DeclareUnicodeCharacter{2218}{\circ}              
\DeclareUnicodeCharacter{221A}{\sqrt}              
\DeclareUnicodeCharacter{221E}{\infty}             
\DeclareUnicodeCharacter{2227}{\land}              
\DeclareUnicodeCharacter{2228}{\lor}               
\DeclareUnicodeCharacter{2229}{\cap}               
\DeclareUnicodeCharacter{2229}{\cap}               
\DeclareUnicodeCharacter{222A}{\cup}               
\DeclareUnicodeCharacter{2243}{\simeq}             
\DeclareUnicodeCharacter{2245}{\cong}              
\DeclareUnicodeCharacter{2260}{\neq}               
\DeclareUnicodeCharacter{2264}{\le}                
\DeclareUnicodeCharacter{2265}{\ge}                
\DeclareUnicodeCharacter{2282}{\subset}            
\DeclareUnicodeCharacter{2283}{\supset}            
\DeclareUnicodeCharacter{2286}{\subseteq}          
\DeclareUnicodeCharacter{2287}{\supseteq}          
\DeclareUnicodeCharacter{2293}{\sqcap}             
\DeclareUnicodeCharacter{2294}{\sqcup}             
\DeclareUnicodeCharacter{22A2}{\vdash}             
\DeclareUnicodeCharacter{22A3}{\dashv}             
\DeclareUnicodeCharacter{22A4}{\top}               
\DeclareUnicodeCharacter{22A5}{\bot}               
\DeclareUnicodeCharacter{22C0}{\bigwedge}          
\DeclareUnicodeCharacter{22C1}{\bigvee}            
\DeclareUnicodeCharacter{22C2}{\bigcap}            
\DeclareUnicodeCharacter{22C3}{\bigcup}            
\DeclareUnicodeCharacter{266D}{\flat}              
\DeclareUnicodeCharacter{266E}{\natural}           
\DeclareUnicodeCharacter{266F}{\sharp}             
\DeclareUnicodeCharacter{2806}{{:}}                
\DeclareUnicodeCharacter{27E6}{\llbracket}         
\DeclareUnicodeCharacter{27E7}{\rrbracket}         
\DeclareUnicodeCharacter{223C}{\sim}               
\DeclareUnicodeCharacter{2248}{\approx}            
\DeclareUnicodeCharacter{2261}{\equiv}             
\DeclareUnicodeCharacter{22C4}{\lozenge}           
\DeclareUnicodeCharacter{22C5}{\Box}               
\DeclareUnicodeCharacter{25FB}{\Box}               
\DeclareUnicodeCharacter{25A1}{\Box}               
\DeclareUnicodeCharacter{2026}{\ldots}             
\DeclareUnicodeCharacter{222B}{\int}               
\DeclareUnicodeCharacter{21C0}{\rightharpoonup}    
\DeclareUnicodeCharacter{2329}{\langle}            
\DeclareUnicodeCharacter{232A}{\rangle}            
\DeclareUnicodeCharacter{2299}{\odot}              
\DeclareUnicodeCharacter{2296}{\ominus}            
\DeclareUnicodeCharacter{2295}{\oplus}             
\DeclareUnicodeCharacter{2298}{\oslash}            
\DeclareUnicodeCharacter{2297}{\otimes}            
\DeclareUnicodeCharacter{2581}{\_}                 
\DeclareUnicodeCharacter{1D41B}{\mathbf}           
\DeclareUnicodeCharacter{1D422}{\textsl}           
\DeclareUnicodeCharacter{1D42B}{\mathrm}           
\DeclareUnicodeCharacter{1D42D}{\text}             
\DeclareUnicodeCharacter{1D42C}{\mathsf}           

\fi
\makeatletter
\newcount\m \newcount\n
\def\twodigits#1{\ifnum #1<10 0\fi \number#1}
\def\hours{\n=\time \divide\n 60
  \m=-\n \multiply\m 60 \advance\m \time
  \twodigits\n:\twodigits\m}
\def\draftfooter{\jobname{} \today{} \hours}
\def\footertext{\draftfooter}
\def\ps@headings{%
  \def\@oddfoot{\hfil \footertext \hfil}
  \let\@evenfoot\@oddfoot
  \def\@oddhead{{\slshape\rightmark}\hfil\thepage}%
  \let\@mkboth\markboth
  \def\chaptermark##1{%
    \markright {\MakeUppercase{%
      \ifnum \c@secnumdepth >\m@ne
        \if@mainmatter
          \@chapapp\ \thechapter. \ %
        \fi
      \fi
      ##1}}}}
\ps@headings
\makeatother
\def\shorttoday{%
  \number\year%
  \ifcase\month\or jan\or feb\or mar\or apr\or may\or
     jun\or jul\or aug\or sep\or oct\or nov\or dec\fi
  \twodigits\day}

\def\sm #1{\begin{smallmatrix}#1\end{smallmatrix}}

\def\bmat#1{\left [\begin{matrix}#1\end{matrix}\right ]}

\def\sm #1{\begin       {smallmatrix}#1\end{smallmatrix}}

\def\beginpicture (#1,#2)(#3,#4){\expr{beginpicture(v(#1,#2),v(#3,#4))}}
\def\beginpictureb(#1,#2)(#3,#4)#5{\expr{beginpicture(v(#1,#2),v(#3,#4),#5)}}

\def\pictaxes{\expr{pictaxes()}}
\def\pictaxes{{\linethickness{0.6pt}\expr{pictaxes()}}}
\def\pictaxes{{\linethickness{0.5pt}\expr{pictaxes()}}}
\def\pictgrid{{\color{GrayPale}\expr{pictgrid()}}}
\def\pictgrid{{\color{GrayPale}\linethickness{0.3pt}\expr{pictgrid()}}}
\def\pictdots#1{\expr{pictdots("#1")}}
\def\pictpiecewise#1{\expr{pictpiecewise("#1")}}
\def\pictpiecewise#1{{\linethickness{1pt}\expr{pictpiecewise("#1")}}}
\def\picturedots(#1,#2)(#3,#4)#5{%
  \vcenter{\hbox{%
  \beginpicture(#1,#2)(#3,#4)%
  \pictaxes%
  \pictdots{#5}%
  \end{picture}%
  }}%
}
\def\picturepiecewise(#1,#2)(#3,#4)#5{%
  \vcenter{\hbox{%
  \beginpicture(#1,#2)(#3,#4)%
  \pictgrid%
  \pictaxes%
  \pictpiecewise{#5}%
  \end{picture}%
  }}%
}

\def\Vector(#1,#2)(#3,#4){\expr{pict2evector(#1, #2, #3, #4)}}

\def\verteq{
 \vcenter{\hbox{%
 \unitlength=1pt%
 \linethickness{0.4pt}%
 \beginpicture(0,0)(5,6)
   \roundcap
   \Line(1,0)(1,6)
   \Line(3,0)(3,6)
 \end{picture}%
 }}%
}

\def\beginpicture(#1,#2)(#3,#4){\expr{beginpicture(v(#1,#2),v(#3,#4))}}
\def\pictaxes{\expr{pictaxes()}}
\def\pictdots#1{\expr{pictdots("#1")}}
\def\picturedotsa(#1,#2)(#3,#4)#5{%
  \vcenter{\hbox{%
  \beginpicture(#1,#2)(#3,#4)%
  \pictaxes%
  \pictdots{#5}%
  \end{picture}%
  }}%
}
\def\picturedots(#1,#2)(#3,#4)#5{%
  \vcenter{\hbox{%
  \beginpicture(#1,#2)(#3,#4)%
  \pictdots{#5}%
  \end{picture}%
  }}%
}

\def\myresizebox#1{%
  \noindent\hbox to \textwidth{\hss
    \resizebox{1.0\textwidth}{!}{#1}%
    \hss}
  }

\def\squigbijtriangle(#1,#2)#3{\polygon*(#1,0)(#2,#3)(#2,-#3)}

\def\defzha#1#2{\expandafter\def\csname zha-#1\endcsname{#2}}
\def\ifzhaundefined#1{\expandafter\ifx\csname zha-#1\endcsname\relax}
\def\zha#1{\ifzhaundefined{#1}
    \errmessage{UNDEFINED ZHA: #1}
  \else
    \csname zha-#1\endcsname
  \fi
}
\def\deftcg#1#2{\expandafter\def\csname tcg-#1\endcsname{#2}}
\def\iftcgundefined#1{\expandafter\ifx\csname tcg-#1\endcsname\relax}
\def\tcg#1{\iftcgundefined{#1}
    \errmessage{UNDEFINED TCG: #1}
  \else
    \csname tcg-#1\endcsname
  \fi
}

\def\defpido#1#2{\expandafter\def\csname pido-#1\endcsname{#2}}
\def\ifpidoundefined#1{\expandafter\ifx\csname pido-#1\endcsname\relax}
\def\pido#1{\ifpidoundefined{#1}
    \errmessage{UNDEFINED PIDO: #1}
  \else
    \csname pido-#1\endcsname
  \fi
}

\def\defpicturedots #1(#2,#3)(#4,#5)#6{%
    \directlua{defpictdots(nil, "#1", #2,#3, #4,#5,nil, "#6")}
  }
\def\defpicturedotsa#1(#2,#3)(#4,#5)#6{%
    \directlua{defpictdots("axes", "#1", #2,#3, #4,#5,nil, "#6")}
  }

\def\picturedotsadef#1(#2,#3)(#4,#5)#6{
  \directlua{ defpictdots("axes", "#1", #2,#3, #4,#5, nil, "#6") }
  \pido{#1}
  }
\def\picturedotsdef #1(#2,#3)(#4,#5)#6{
  \directlua{ defpictdots(nil,    "#1", #2,#3, #4,#5, nil, "#6") }
  \pido{#1}
  }

\def\defub#1#2{\expandafter\def\csname ub-#1\endcsname{#2}}
\def\ifubundefined#1{\expandafter\ifx\csname ub-#1\endcsname\relax}
\def\ub#1{\ifubundefined{#1}
    \errmessage{UNDEFINED UB: #1}
  \else
    \csname ub-#1\endcsname
  \fi
}

\usepackage[backend=biber,
   style=alphabetic]{biblatex}            
\begin{document}


\ifluatex
  \catcode`\^^J=10
  \directlua{dofile "dednat6load.lua"}
\else
  \input\jobname.dnt
  \def\pu{}
\fi

\def\co#1{{%
  \def\%{\char37}%
  \def\\{\char92}%
  \def\^{\char94}%
  \def\~{\char126}%
  \tt#1%
  }}
\def\qco#1{`\co{#1}'}
\def\qqco#1{``\co{#1}''}

\def\respcomp{\mathsf{respcomp}}
\def\respids {\mathsf{respids}}
\def\sqcond  {\mathsf{sqcond}}
\def\assoc   {\mathsf{assoc}}
\def\idL     {\mathsf{idL}}
\def\idR     {\mathsf{idR}}
\def\univ    {\mathsf{univ}}
\def\Ran     {\mathsf{Ran}}

\def\sfC  {\mathsf{C}}
\def\sfSet{\mathsf{Set}}
\def\Ring {\mathbf{Ring}}
\def\nameof#1{\ulcorner#1\urcorner}
\def\catK {\mathbf{K}}

\def\veq{\rotatebox{90}{$=$}}
\def\liml{\underleftarrow {\lim}{}}
\def\limr{\underrightarrow{\lim}{}}

\def\DONE{(DONE)}
\def\DONE{}

%

\long\def\ColorRed   #1{{\color{Red1}#1}}
\long\def\ColorViolet#1{{\color{MagentaVioletLight}#1}}
\long\def\ColorViolet#1{{\color{Violet!50!black}#1}}
\long\def\ColorGreen #1{{\color{SpringDarkHard}#1}}
\long\def\ColorGreen #1{{\color{SpringGreen4}#1}}
\long\def\ColorGreen #1{{\color{SpringGreenDark}#1}}
\long\def\ColorGray  #1{{\color{GrayLight}#1}}
\long\def\ColorGray  #1{{\color{black!30!white}#1}}
\long\def\ColorBrown #1{{\color{Brown}#1}}
\long\def\ColorBrown #1{{\color{brown}#1}}

%

\title{On my favorite conventions for drawing the missing diagrams in
  Category Theory}

\author{Eduardo Ochs}


\maketitle


\begin{abstract}

  I used to believe that my conventions for drawing diagrams for
  categorical statements could be written down in one page or less,
  and that the only tricky part was the technique for reconstructing
  objects ``from their names'' (sec.\ref{to-deserve-a-name})... but
  then I found out that this is not so.

  This is an attempt to explain, with motivations and examples, all
  the conventions behind a certain diagram, called the ``Basic
  Example'' in the text. Once the conventions are understood that
  diagram becomes a ``skeleton'' for a certain lemma related to the
  Yoneda Lemma, in the sense that both the statement and the proof of
  that lemma can be reconstructed from the diagram. The last sections
  discuss some simple ways to extend the conventions; we see how to
  express in diagrams the (``real'') Yoneda Lemma and a corollary of
  it, how to define comma categories,
  %
  %
  and how to formalize the diagram for ``geometric morphism for
  children'' mentioned in sec.\ref{missing-diagrams}.

  \msk

  People in CT usually only share their ways of visualizing things
  when their diagrams cross some threshold of of mathematical
  relevance --- and this usually happens when they prove new theorems
  with their diagrams, or when they can show that their diagrams can
  translate calculations that used to be huge into things that are
  much easier to visualize. The diagrammatic language that I present
  here lies below that threshold --- and so it is a ``private''
  diagrammatic language, that I am making public as an attempt to
  establish a dialogue with other people who have also created their
  own private diagrammatic languages.

\end{abstract}

\newpage

%

\renewcommand{\cfttoctitlefont}{\bfseries}
\setlength{\cftbeforesecskip}{2.5pt}


%

\tableofcontents

%

\section{Missing diagrams}
\label{missing-diagrams}

I need to tell a long story here.

Let me start with some quotes. This one is from Eilenberg and Steenrod
(\cite[p.ix]{EilenbergSteenrod}, but I learned it from
\cite[pp.82--83]{Kromer}):

\begin{quotation}

  The diagrams incorporate a large amount of information. Their use
  provides extensive savings in space and in mental effort. In the
  case of many theorems, the setting up of the correct diagram is the
  major part of the proof. We therefore urge that the reader stop at
  the end of each theorem and attempt to construct for himself the
  relevant diagram before examining the one which is given in the
  text. Once this is done, the subsequent demonstration can be
  followed more readily; in fact, the reader can usually supply it
  himself.

\end{quotation}

I spent a {\sl lot} of my time studying Category Theory trying to
``supply the diagrams myself''. In \cite{EilenbergSteenrod} supplying
the diagrams is not very hard (I guess), but in books like
\cite{CWM2}, in which most important concepts involve several
categories, I had to rearrange my diagrams hundreds of times until I
reached ``good'' diagrams...


The problem is that I expected too much from ``good'' diagrams. The
next quotes are from the sections 1 and 12 of an article that I wrote
about that (\cite{IDARCT}):

\begin{quotation}

  My memory is limited, and not very dependable: I often have to
  rededuce results to be sure of them, and I have to make them fit in
  as little ``mental space'' as possible...

  Different people have different measures for ``mental space'';
  someone with a good algebraic memory may feel that an expression
  like $\mathsf{Frob}: Σ_f(P ∧ f^* Q) ≅ Σ_f P ∧ Q$ is easy to
  remember, while I always think diagramatically, and so what I do is
  that I remember this diagram,

  \begin{center}
  \includegraphics[height=60pt]{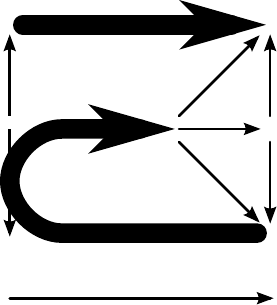}
  \end{center}

  \noindent and I reconstruct the formula from it.

\end{quotation}

\begin{quotation}

  Let's call the ``projected'' version of a mathematical object its
  ``skeleton''. The underlying idea in this paper is that for the
  right kinds of projections, and for some kinds of mathetical
  objects, it should be possible to reconstruct enough of the original
  object from its skeleton and few extra clues --- just like
  paleontologists can reconstruct from a fossil skeleton the look of
  an animal when it was alive.

\end{quotation}

I was searching for a diagrammatic language that would let me express
the ``skeletons'' of categorical definitions and proofs. I wanted
these skeletons to be easy to remember --- partly because they would
have shapes that were easy to remember, and partly because they would
be similar to ``archetypal cases'' (\cite[section 16]{IDARCT}).

\bsk

In 2016 and 2017 I taught a seminar course for undergraduates that
covered a bit of Category Theory in the end --- see Section
\ref{teaching-adjunctions} and \cite{OchsWLD2019} --- and this forced
me to invent new techniques for working in two different styles in
parallel: a style ``for adults'', more general, abstract, and formal,
and another ``for children'', with more diagrams and examples. After
some semesters, and after writing most of the material that became
\cite{PH1}, I tried to read again some parts of Johnstone's ``Sketches
of an Elephant'', a book that always felt quite impenetrable to me,
and I found a way to present geometric morphisms in toposes to
``children''. It was based on this diagram,
\pu
$$
  \diag{gm-for-adults}
  \qquad
  \def\LG{\pshAargs{G_2}{G_3}{G_4}{G_5}}
  \def\G {\pshBargs{G_1}{G_2}{G_3}{G_4}{G_5}{G_6}}
  \def\H {\pshAargs{H_2}{H_3}{H_4}{H_5}}
  \def\RH{\pshBargs{H_2{×_{H_4}}H_3}{H_2}{H_3}{H_4}{H_5}{1}}
  \scalebox{0.6}{$
  \diag{gm-for-children}
  $}
$$
that we will discuss in detail in \ref{gms-for-children}. Its left
half is a generic geometric morphism (``for adults''), and its right
half is a very specific geometric morphism (``for children'') in which
everything is easy to understand and to visualize, and that turns out
to be ``archetypal enough''.

I showed that to the few categorists with whom I had contact and the
feedback that I got was quite positive. A few of them --- the ones who
were strictly ``adults'' --- couldn't understand why I was playing
with particular cases, and even worse, with finite categories, instead
of proving things in the most general case possible, but some others
said that these ideas were very nice, that they knew a few bits about
geometric morphisms but those bits didn't connect well, and that now
they had a family of particular cases to think about, and they had
much more intuition than before.

That was the first time that my way of using diagrams yielded
something so nice! This was the excuse that I needed to organize a
workshop on diagrammatic languages and ways to use particular
cases; here's how I advertised it (from \cite{OchsLucatelli}):
\begin{quotation}

  When we explain a theorem to children --- in the strict sense of the
  term --- we focus on concrete examples, and we avoid
  generalizations, abstract structures and infinite objects.

  When we present something to ``children'', in a wider sense of the
  term that means ``people without mathematical maturity'', or even
  ``people without expertise in a certain area'', we usually do
  something similar: we start from a few motivating examples, and then
  we generalize.

  One of the aims of this workshop is to discuss techniques for {\sl
    particularization} and {\sl generalization}. Particularization is
  easy; substituing variables in a general statement is often enough
  to do the job. Generalization is much harder, and one way to
  visualize how it works is to regard particularization as a
  projection: a coil projects a circle-like shadow on the ground, and
  we can ask for ways to ``lift'' pieces of that circle to the coil
  continously. {\sl Projections} lose dimensions and may collapse
  things that were originally different; {\sl liftings} try to
  reconstruct the missing information in a sensible way. There may be
  several different liftings for a certain part of the circle, or
  none. Finding good generalizations is somehow like finding good
  liftings.

  The second of our aims is to discuss {\sl diagrams}. For example, in
  Category Theory statements, definitions and proofs can be often
  expressed as diagrams, and if we start with a general diagram and
  particularize it we get a second diagram with the same shape as the
  first one, and that second diagram can be used as a version ``for
  children'' of the general statement and proof. Diagrams were for a
  long time considered second-class entities in CT literature
  (\cite{Kromer} discusses some of the reasons), and were omitted;
  readers who think very visually would feel that part of the work
  involved in understanding CT papers and books would be to
  reconstruct the ``missing'' diagrams from algebraic statements.
  Particular cases, even when they were the motivation for the general
  definition, are also treated as somewhat second-class --- and this
  inspires a possible meaning for what can call ``Category Theory for
  Children'': to start from the diagrams for particular cases, and
  then ``lift'' them to the general case. Note that this can be done
  outside Category Theory too; \cite{Jamnik} is a good example.

  Our third aim is to discuss {\sl models}. A standard example is that
  every topological space is a Heyting Algebra, and so a model for
  Intuitionistic Predicate Logic, and this lets us explain visually
  some features of IPL. Something similar can be done for some modal
  and paraconsistent logics; we believe that the figures for that
  should be considered more important, and be more well-known.

\end{quotation}

This is from the second announcement:

\begin{quotation}

  If we say that categorical definitions are ``for adults'' - because
  they may be very abstract - and that particular cases, diagrams, and
  analogies are ``for children'', then our intent with this workshop
  becomes easy to state. ``Children'' are willing to use ``tools for
  children'' to do mathematics, even if they will have to translate
  everything to a language ``for adults'' to make their results
  dependable and publishable, and even if the bridge between their
  tools ``for children'' and ``for adults'' is somewhat defective,
  i.e., if the translation only works on simple cases...

  We are interested in that {\sl bridge} between maths ``for adults''
  and ``for children'' in several areas. Maths ``for children'' are
  hard to publish, even informally as notes (see this thread

  \msk

  \centerline{\url{http://angg.twu.net/categories-2017may02.html}}

  \msk

  \noindent in the Categories mailing list), so often techniques are
  rediscovered over and over, but kept restricted to the ``oral
  culture'' of the area.

  Our main intents with this workshop are:

  \begin{itemize}

     \item to discuss (over coffe breaks!) the techniques of the
       ``bridge'' that we currently use in seemingly ad-hoc ways,

     \item to systematize and ``mechanize'' these techniques to make
       them quicker to apply,

     \item to find ways to publish those techniques --- in journals or
       elsewhere,

     \item to connect people in several areas working in related
       ideas, and to create repositories of online resources.

  \end{itemize}

\end{quotation}

In the UniLog 2018 I was able to chat with several categorists, and
they told me about the oral culture of CT and showed me that it was
not as I was guessing, and I also spent two evenings with Peter Arndt
working on factorizations of geometric morphisms ``for children'' ---
and this made me feel that I could present applications of this
diagrammatic language in conferences that were more top-level-ish in
some sense.

The following quote is from the abstract of my submission (\cite{MDE})
to the ACT2019:
\begin{quotation}

  Imagine two category theorists, Aleks and Bob, who both think very
  visually and who have exactly the same background. One day Aleks
  discovers a theorem, $T_1$, and sends an e-mail, $E_1$, to Bob,
  stating and proving $T_1$ in a purely algebraic way; then Bob is
  able to reconstruct by himself Aleks's diagrams for $T_1$ exactly as
  Aleks has thought them. We say that Bob has reconstructed the
  {\it missing diagrams} in Aleks's e-mail.

  Now suppose that Carol has published a paper, $P_2$, with a theorem
  $T_2$. Aleks and Bob both read her paper independently, and both
  pretend that she thinks diagrammatically in the same way as them.
  They both ``reconstruct the missing diagrams'' in $P_2$ in the same
  way, even though Carol has never used those diagrams herself.

\end{quotation}
and this from my submission (\cite{OchsTallinnAbs}) to Diagrams 2020:
\begin{quotation}

  Category Theory gives the impression of being an area where most
  concepts and arguments are stated and formalized via diagrams, but
  this is not exactly true... in most texts almost everything is done
  algebraically, and the reader is expected to be able to reconstruct
  the ``missing diagrams'' by himself.

  I used to believe, as an outsider, that some people who grew up
  immersed the oral culture of the area would know several techniques
  for ``drawing the missing diagrams''. My main intent when I
  organized the workshop ``Logic for Children'' at the UniLog 2018
  \cite{OchsLucatelli} was to collect some of these folklore
  techniques, compare them with the ones that I had developed myself
  to study CT, and formalize them all --- but what I found instead was
  that everybody that I could get in touch with used their own ad-hoc
  techniques, and that what I was trying to do was either totally new
  to them, or at least new in its level of detail.

\end{quotation}

The story continues in the last three sections --- that also explains
why I decided to write these notes using the first person in most
places.




%

\section{The conventions \DONE}
\label{the-conventions}

The conventions that I will present now are the ones that we need for
this diagram (called the ``Basic Example'' from here on), that is
essentially the Proposition 1 in the proof of the Yoneda Lemma in
\cite[Section III.2]{CWM2}:
%
%
\pu
$$\scalebox{2.0}{$
  \diag{Basic-Example}
  $}
$$

\begin{itemize}

\item[(CD)] Our diagrams are made of components that are nodes and
  arrows. The nodes can contain arbitrary expressions. The arrows work
  as connectives, and each arrow can be interpreted as the top-level
  connective in the smallest subexpression that contains it. For
  example, the curved arrow in the diagram above can be interpreted
  as:
  $$(A\ton{η}RC)↔(\catB(C,-) \ton{T} \catA(A,R-)).
  $$

\item[(C$→$)] Arrows that look like `$→$' (\qqco{\\to}) represent
  hom-sets, or, in $\Set$, spaces of functions. When a `$→$' arrow is
  named the name stands for an element of that hom-set. For example,
  in $A \ton{η} RC$ we have $η:A→RC$.

\item[(C$↦$)] Arrows that look like `$↦$' (\qqco{\\mapsto}) represent
  internal views of functions or functors. This has some subtleties;
  see Section \ref{internal-views}.

\item[(C$↔$)] Arrows that look like `$↔$' (\qqco{\\leftrightarrow})
  represent bijections or isomorphisms.

\item[(CAI)] ``Above'' usually means ``inside'', or ``internal view''.
  In the diagram above the morphism $η:A→RC$ is in $\catA$ and $C$ is
  an object of $\catB$. Also, the arrow $C \mapsto RC$ is above $\catB
  \ton{R} \catA$, and this means that it is an internal view of the
  functor $R$. Note that {\sl usually} is not {\sl always} --- and
  $\catB \ton{R} \catA$ is not an internal view of $\catB(C,-) \ton{T}
  \catA(A,R-)$.

\item[(CO)] When the definition of a component of our diagram is
  ``obvious'' in the sense of ``there is a unique natural construction
  for an object with that name'', we will usually omit its definition
  and {\sl pretend} that it is obvious; same for its uniqueness. See
  Section \ref{to-deserve-a-name}.

\item[(CC)] Everything commutes by default, and non-commutative cells
  have to be indicated explicitly. See Section \ref{freyd-notation}.

\item[(CTL)] The default ``meaning'' for a diagram is the definition
  of its top-level component. There is a natural partial order on the
  components of a diagram, in which $α \prec β$ iff $α$ is ``more
  basic'' than $β$, or, in other words, if $α$ needs to be defined
  before $β$. In the diagram above the top-level component is the
  curved bijection.

\item[(CAdj)] {\sl I use shapes based on my way of drawing adjunctions
  whenever possible.} I like adjunctions so much that when I want to
  explain Category Theory to someone who knows just a little bit of
  Maths I always start by the adjunction $({×}B)⊣(B{→})$ of Section
  \ref{internal-view-adjunction}; I always draw it in a canonical way,
  with the left adjoint going left, the right adjoint going right, and
  the morphisms going down. In Proposition 1 of \cite[Section
    III.2]{CWM2} the map $η$ is a universal arrow, and someone who
  learns adjunctions first sees the unit maps $η:A→(B{→}(A{×}B))$ as
  the first examples of universal arrows --- so that's why the upper
  part of the diagram above is drawn in this position.

\item[(COT)] We use a notation as close to the original text as
  possible, especially when we are trying to draw the missing diagrams
  for some existing text. If we were drawing the missing diagrams for
  the Proposition 1 of \cite[Section III.2]{CWM2} our diagram would be
  this:
  %
%
$$\pu
  \diag{yoneda-cwm-0-small}
$$
but I hate Mac Lane's choice of letters, so I decided to use another
notation here.

\item[(CSk)] Suppose that we have a piece of text --- say, a paragraph
  $P$ --- and we want to reconstruct the ``missing diagram'' $D$ for
  $P$. Ideally this $D$ should be a ``skeleton'' for $P$, in the sense
  that it should be possible to reconstruct the ideas in $P$ from the
  diagram $D$ using very few extra hints; see \cite[sec.12]{IDARCT}.

\item[(CFSh)] The image by a functor of a diagram $D$ is drawn with
  the same shape as $D$.

\item[(CISh)] The internal view of a diagram $D$ is drawn with the
  same shape as $D$, modulo duplications --- see section
  \ref{internal-views}.

\item[(CPSh)] A particular case of a diagram $D$ is drawn with the
  same shape as $D$.

\item[(CNSh)] A translation of a diagram $D$ to another notation is
  drawn with the same shape as $D$.

\end{itemize}


Note that I have presented these conventions in a human-friendly way,
that is somewhat informal and admits exceptions and extensions. Some
simple examples of extensions will be discussed in Section
\ref{extensions}.

See \cite{PenroseSIGGRAPH2020} for a system that produces diagrams
from conventions and specifications and then lets the user adjust
these generated diagrams to make them clearer and more aesthetically
pleasing --- but as far as I know Penrose can only {\sl produce}
diagrams, not {\sl read} them.

%
\section{Finding ``the'' object with a given name \DONE}
\label{to-deserve-a-name}

One of the books that I tried to read when I was starting to learn
Category Theory was Mac Lane's \cite{CWM2}. It is written for readers
who know a lot of mathematics and who can follow some steps that it
treats as obvious. I was not (yet) a reader like that, but I wanted to
become one.

There is one specific thing that \cite{CWM2} does pretending that it
is obvious that I found especially fascinating. It ``defines''
functors by describing their actions on objects, and it leaves to the
reader the task of discovering their actions on morphisms. Let's see
how to find these actions on morphisms.

A functor $F:\catA→\catB$ has four components:
$$F=(F_0, F_1, \respids_F, \respcomp_F).$$
They are its action on objects, its action on morphisms, the assurance
that it takes identity maps to identity maps, and the assurance that
it respects compositions. When Mac Lane says this,
\begin{quote}
Fix a set $B$. Let $(×B)$ denote {\sl the} functor that takes each set
$A$ to $A×B$.
\end{quote}
he is saying that $(×B)_0 A = A×B$, or, more precisely, this:
$$(×B)_0 := λA.\,A×B$$

The ``{\sl the}'' in the expression ``Let $(×B)$ denote {\sl the}
functor...'' implies that the precise meaning of $(×B)_1$ is easy to
find, and that it is easy to prove $\respids_{(×B)}$ and
$\respcomp_{(×B)}$.

If $f:A'→A$ then $(×B)_1 f : (×B)_0 A' → (×B)_0 A$. We know the {\sl
  name} of the image morphism, $(×B)_1 f$, and its {\sl type},
$$(×B)_1 f : A'×B → A×B,$$
and it is implicit that there is an ``obvious'' natural construction
for this $(×B)_1 f$ from $f$. A natural construction is ---
TA-DAAAA!!! --- a $λ$-term, so we are looking for a term of type $A'×B
→ A×B$ that can be constructed from $f:A'→A$.

In a big diagram:
\pu
$$\ded{foo1} \quad ⇒ \ded{foo2}$$

A double bar in a derivation means ``there are several omitted steps
here'', and sometimes a double bar suggests that these omitted steps
are obvious. The derivation on the left says that there is an
``obvious'' way to build a $(×B)_1f:A'{×}B→A{×}B$ from a
``hypothesis'' $f:A'→A$. If we expand its double bar we get the tree
at the right, that shows that the ``precise meaning'' for $(×B)_1f$ is
$(λp⠆A'{×}B.(f(πp),π'p)$. More formally (and erasing a typing),
$$(×B)_1 := λf.(λp.(f(πp),π'p)).$$

The expansion of the double bar above becomes something more familiar
if we translate the trees to Logic using Curry-Howard:
\pu
$$\ded{foo-logic1} \quad ⇒ \ded{foo-logic2}$$

We obtain the tree at the right by {\sl proof search}.

Let's give a name for the operation above that obtained a term of type
$A'×B→A×B$: we will call that operation {\sl term search}, or, as it
is somewhat related to type inference, {\sl term inference}.

Term search may yield several different construction and trees, and so
several non-equivalent terms of the desired type. When Mac Lane says
``{\sl the} functor $(×B)$'' he is indicating that:

\begin{itemize}

\item a term for $(×B)_1$ is easy to find (note that we use the
  expression ``a {\sl precise meaning} for $(×B)_1$''),

\item all other natural constructions for something that ``deserves
  the name'' $(×B)_1$ yield terms equivalent to that first, most
  obvious one,

\item proving $\respids_{(×B)}$ and $\respcomp_{(×B)}$ is trivial.

\end{itemize}

In many situations we will start by just the name of a functor, as the
``$(×B)$'' in the example above, and from that name it will be easy to
find {\sl the} ``precise meaning'' for $(×B)_0$, and from that the
``precise meaning'' for $(×B)_1$, and after that proofs that
$\respids_{(×B)}$ and $\respcomp_{(×B)}$. We will use the expression
``...deserving the name...'' in this process --- terms for $(×B)_0$,
$(×B)_1$, $\respids_{(×B)}$, and $\respcomp_{(×B)}$ ``deserve their
names'' if they obey the expected constraints.

For a more thorough discussion see \cite{IDARCT}.

\msk

{\sl Note: I am not aware of any papers or books that discuss how to
  (re)construct a functor from its action on objects, or from its
  name. If you have any references, please let me know!}

\msk

These ideas of ``finding a precise meaning'' and ``finding (something)
deserving that name'' can also be applied to morphisms, natural
transformations, isomorphisms, and so on.

In Section \ref{basic-example-bij} we will see how to find natural
constructions for the two directions of the bijection in the Basic
Example --- or how the expand the double bars in the two derivations
here:
$$\pu
  \diag{Basic-Example}
  \qquad
  \begin{array}{c}
  \ded{bij-down}
  \\ \\
  \ded{bij-up}
  \end{array}
$$


%

\section{Freyd's diagrammatic language \DONE}
\label{freyd-notation}

In \cite{Freyd76} Peter Freyd presents a very nice diagrammatic
language that can be used to express {\sl some} definitions from
Category Theory. For example, this is the statement that a category
has all equalizers:
$$\pu
  \scalebox{0.8}{$
  \diag{cat-has-equalizers}
  $}
$$

All cells in these diagrams commute by default, and non-commuting
cells have to be indicated with a `?'. Each vertical bar with a `$∀$'
above it means ``for all extensions of the previous diagram to this
one such that everything commutes''; a vertical bar with a `$∃!$'
above it means ``there exists a unique extension of the previous
diagram to this one such that everything commutes'', and so on. See
the scan in \cite{Freyd76} for the basic details of how to formalize
these diagrams, and the book \cite[p.28 onwards]{FreydScedrov}, for
tons of extra details, examples, and applications.

Let's call the subdiagrams of a diagram like the one above its
``stages''. Its stage 0 is empty, its stage 1 has two objects and two
arrows, its last stage has four objects and five arrows, and the
quantifiers separating the stages are $Q_1=∀$, $Q_2=∃$, $Q_3=∀$,
$Q_4=∃!$. They are structured like this:
%
\pu
$$
  \diag{freyd-stages-1}
$$

I was not very good at drawing all stages separately --- it was
boring, it took me too long, and I often got distracted and committed
errors --- so I started to play with extensions of that diagrammatic
language.

%
\subsection{Adding quantifiers \DONE}
\label{freyd-with-quantifiers}

Here is a simple way to draw all stages at once. We start from a
diagram for the ``last stage with quantifiers'', that we will call
$LSQ$:
%
$$\pu
  \scalebox{1.75}{$
  \diag{cat-has-equalizers-with-quants}
  $}
$$

We can recover all the stages and quantifiers from it. The numbered
quantifiers in it are $∀_1$, $∃_2$, $∀_3$, and $∃!_4$. The highest
number in them 4, so we set $n=4$ ($n$ is the index of the last
stage), and we set ``stage 4 with quantifiers'', $SQ_4$, to $LSQ$. To
obtain the $SQ_3$ from $SQ_4$ we delete all nodes an arrows in $SQ_4$
that are annotated with a `$∃!_4$'; to obtain $SQ_2$ from $SQ_3$ we
delete all nodes an arrows in $SQ_3$ that are annotated with a
`$∀_3$', and so on until we get a diagram $SQ_0$, that in this example
is empty. To obtain each $S_k$ --- a stage in the original
diagrammatic language from Freyd, that doesn't have quantifiers ---
from the corresponding $SQ_k$ we treat all the quantifiers in $SQ_k$
as mere annotations, and we erase them; for example, `$∃_2e$' becomes
`$e$', and $∀_1A$ becomes $A$. To obtain the quantifiers $Q_1$, $Q_2$,
$Q_3$, $Q_4$ that are put in the vartical bars that separate the
stages, we just assign $∀_1$, $∃_2$, $∀_3$, and $∃!_4$ to them,
without the numbers in the subscripts.

Bonus convention: when the quantifiers in a diagram are just `$∀$'s
and `$∃!$'s without subscripts the `$∀$'s are to be interpreted as
`$∀_1$' and the `$∃!$'s as `$∃!_2$'s.

%
\subsection{Adding functors \DONE}
\label{freyd-with-functors}

Freyd's language can't represent functors\footnote{As far as I know
  --- I don't know \cite{FreydScedrov} very well.}, and I wanted to
use it to draw the missing diagrams for definitions involving
functors, so I had to extend it again.

Let me use an example to discuss this. This is the definition of
universal arrow in \cite[p.55]{CWM2}, including the original diagram,
modulo change of letters:

%

\begin{quotation}

  {\bf Definition.} If $R: \catB→\catA$ is a functor and $A$ an object
  of $\catA$, a universal arrow from $A$ to $R$ is a pair $(B,η)$
  consisting of an object $B$ of $\catB$ and and arrow $η:A→RB$ of
  $\catA$ such that to every pair $(B',g)$ with $B'$ an object of
  $\catB$ and $g:A→RB'$ an arrow of $\catA$, there is a unique arrow
  $f:B→B'$ of $\catB$ with $Rf∘η=g$. In other words, every arrow $h$
  to $R$ factors uniquely through the universal arrow $η$, as in the
  commutative diagram:
  $$\pu
    \diag{univ-arrow-cwm-my-letters}
  $$

\end{quotation}

The definition itself goes only up to the ``with $Rf∘η=g$.'', so let
me ignore the part starting from ``In other words'', and draw a better
``missing diagram'' for the definition:
%
$$\pu
  \diag{universal-arrow-stages}
$$

This diagram is quite close to being a skeleton for the definition of
universal arrow. It can be interpreted as a proposition, and the only
extra hint that we need is that ``universalness'' for the arrow $η$
corresponds to the truth of that proposition. Here's how to extract
the proposition from it:
$$\begin{array}{rl}
  \text{In a context where:}
    & \catA \text{ is a category}, \\
    & \catB \text{ is a category}, \\
    & R:\catB \to \catA, \\
    & A ∈ \catA, \\
    & B ∈ \catB, \\
    & η:A→RB, \\
  \text{for all}
    & B'∈\catB \text{ and} \\
    & g:A→RB', \\
  \text{there exists a unique}
    & f:B→B' \text { such that} \\
    & Rf∘η=g. \\
  \end{array}
$$

To convert that to a definition of universalness we just have to
replace the ``for all'' by ``$(B,η)$ is a universal arrow for $A$ to
$R$ iff for all''.

The convention for quantifiers from sec.\ref{freyd-with-quantifiers}
lets us rewrite the diagram in three stages above as:
%
$$\pu
  \scalebox{1.5}{$
  \diag{universal-arrow-quants}
  $}
$$

Also, I noticed that I could omit most typings when they could be
inferred from the diagram. I could ``formalize'' the diagram above as:
``in a context where $(\catA, \catB, R, A, B, η)$ are as in the
diagram above, we say that $(B,η)$ is a universal arrow from $A$ to
$R$ when $∀(B',g).∃!f.(Rf∘η=g)$''. This may be too loaded to be used
in public, but it's very practical for private notes --- and I can
even omit the ``$Rf∘η=g$'', as everything commutes by default.

\bsk

Note that when we erase a node or arrow we also erase everything that
depends on it. In the example above $SQ_2$ has an arrow labeled
$∃!_2f$; to obtain $SQ_1$ from $SQ_2$ we have to erase that arrow, the
arrow $Rf$, and the arrow $f \mapsto Rf$ --- and to obtain $SQ_0$ from
$SQ_1$ we have to erase the arrow $g$, the node $B'$, the node $RB'$,
and the arrow $B' \mapsto RB'$.


%
\section{Internal views \DONE}
\label{internal-views}

My favorite way of introducing internal views is with the diagram
below:
%
%
\def\ooo(#1,#2){\begin{picture}(0,0)\put(0,0){\oval(#1,#2)}\end{picture}}
\def\oooo(#1,#2){{\setlength{\unitlength}{1ex}\ooo(#1,#2)}}
%
\pu
$$\begin{array}{rrcl}
   \sqrt{\;\;}: & \N &→& \R \\
                &  n &↦& \sqrt{n} \\
   \end{array}
   \qquad
   \diag{second-blob-function}
$$

\def\longmapsto{\diagxyto/|->/}

The parts with the two blobs and `$\longmapsto$'s between them is
based on how I was taught sets and functions when I was a kid; it is
an internal view of the $\N \ton{\sqrt{\phantom{a}}} \R$ below it. Not
all elements of $\N$ are shown in the blob-view of $\N$, but the ones
that are shown are named; compare this with \cite[p.2
  onwards]{LawvereRosebrugh}, in which the elements are usually dots.

The arrow $n \longmapsto \sqrt{n}$ between the blobs shows a {\sl
  generic element} of $\N$ and its image, and the other
`$\longmapsto$'s are {\sl substitution instances of it}, like this:
$$(n \longmapsto \sqrt{n}) [n:=2] = (2 \longmapsto \sqrt{2})$$

In some cases, like $4 \longmapsto 2$, we write 2 instead of
$\sqrt{4}$ because $\sqrt{4}$ ``reduces to'' 2, as explained in the
next section.

%
\subsection{Reductions \DONE}
\label{reductions}

\def\squigton#1{\overset{#1}{\squigto}}

The convention (C$\mapsto$) says that an arrow $α \mapsto β$ above an
arrow $A \ton{f} B$ should be interpreted as meaning $f(α) \squigto
β$, where `$\squigto$' means ``reduces to''; the standard example is
$\sqrt{4} \squigto 2$. In a diagram:
%
$$\pu
  \diag{reductions-0}
$$

The idea of reduction comes from $λ$-calculus. We write $α
\squigton{1} β$ to say that the term $α$ reduces to $β$ in one step,
and $α \squigton{*} γ$ to say that there is a finite sequence of
one-step reductions that reduce $α$ to $γ$. Here we are interested in
reduction in a system with constants, in which for example
$(\sqrt{\phantom{a}})(4) \squigton{1} 2$.

Here is a directed graph that shows all the one-step reductions
starting from $g(2+3)$, considering $g(a) = a·a+4$:
%
%
$$\pu
  \diag{reduce-g}
$$

Note that all reductions sequences starting from $g(2+3)$ terminate at
the same term, 29 --- ``the term $g(2+3)$ is strongly normalizing''
---, and reduction sequences from $g(2+3)$ may ``diverge'' but they
``converge'' later --- this is the ``Church-Rosser Property'', a.k.a.
``confluence''.

A good place to learn about reduction in systems with constants is
\cite{SICP}.



%
\subsection{Functors \DONE}
\label{internal-view-functor}


By the convention (CFSh) the image of the diagram above $\catA$ in the
diagram below --- remember that {\sl above} usually means {\sl inside}
---
%
$$\pu
  \diag{internal-view-functor-0}
$$
is a diagram with the same shape over $\catB$. We draw it like this:
%
$$\pu
  \diag{internal-view-functor-1}
$$

In this case we don't draw the arrows like $A_1 \mapsto FA_1$ because
there would be too many of them --- we leave them implicit.

We say that the diagram above is {\sl an} internal view of the functor
$F$. To draw {\sl the} internal view of the functor $F: \catA → \catB$
we start with a diagram in $\catA$ that is made of two generic objects
and a generic morphism between them. We get this:

$$\pu
  \diag{internal-view-functor-2}
$$

Compare this with the diagram with blob-sets in Section
\ref{internal-views}, in which the `$n \mapsto \sqrt{n}$' says where a
generic element is taken.

Any arrow of the form $α \mapsto β$ above a functor arrow $\catA
\ton{f} \catB$ is interpreted as saying that $F$ takes $α$ to $β$, or,
in the terminology of the section \ref{reductions}, that $Fα$ reduces
to $β$. So this diagram 
%
$$\pu
  \diag{internal-view-functor-3}
$$
defines $(A×)$ as:
$$\begin{array}{rcl}
  (A×)_0 &:=& λB.\,A×B,\\
  (A×)_1 &:=& λf.λp.(πp,f(π'p)).\\
  \end{array}
$$

In this case we can also use internal views of $(A×)$ to define
$(A×)_1$:
%
$$\pu
  \diag{internal-view-functor-4}
$$

%
\subsection{Natural transformations \DONE}
\label{internal-view-NT}


Suppose that we have two functors $F,G:\catA → \catB$ and a natural
transformation $T:F→G$. A first way to draw an internal view of $T$ is
this:
%
$$\pu
  \diag{internal-view-NT-0}
$$

If we start with a morphism $h:C→D$ in $\catA$, like this,
%
%
$$\pu
  \diag{NT-00}
$$
the convention (CFSh) would yield an image of $h$ by $F$ and another
by $G$, and we can draw the arrows $TC$ and $TD$ to obtain a commuting
square in $\catB$:
%
%
$$\pu
  \diag{NT-0}
$$

This way of drawing internal views of natural transformations yields
diagrams that are too heavy, so we will usually draw them as just
this:
%
%
$$\pu
  \diag{NT-1}
$$
Note that the input morphism is at the left, and above $F \ton{T} G$
we draw its images by $F$, $G$, and $T$.

When the codomain of $F$ and $G$ is $\Set$ we will sometimes also draw
at the right an internal view of the commuting square, like this:
%
$$\pu
  \diag{NT-2}
$$
Then the commutativity of the middle square is equivalent to
$∀x∈FC.(Gh∘TC)(x)=(TD∘Ff)(x)$. Note that in this case the square at
the right is an internal view of an internal view.

In Section \ref{to-deserve-a-name} we saw that a functor has four
components. A natural transformation has two: $T=(T_0, \sqcond_T)$,
where $T_0$ is the operation $C \mapsto TC$ and $\sqcond_T$ is the
guarantee that all the induced squares commute. Sometimes we will use
the upper line of the internal view of the internal view to define
$T_0$ --- see Section \ref{basic-example-NT} for an example of this.


%
\subsection{Adjunctions \DONE}
\label{internal-view-adjunction}


We will draw adjunctions like this,
%
$$\pu
  \diag{generic-adjunction}
$$
with the left adjoint going left and the right adjoint going right. My
favorite names for the left and right adjoints are $L$ and $R$. The
standard notation for that adjunction is $L⊣R$.

The top-level component of the diagram above is the bijection arrow in
the middle of the square --- it says that $\Hom(LA,B) ↔ \Hom(A,RB)$.
It is implicit that we have bijections like that for all $A$ and $B$;
it is also implicit that that bijection is natural in some sense.

We will sometimes expand adjunction diagrams by adding unit and counit
maps, the unit and the unit as natural transformations, the actions of
$L$ and $R$ on morphisms, and other things. For example:
%
$$\pu
  \diag{generic-adjunction-with-with}
$$

We can obtain the naturality conditions by regarding $♭$ and $♯$ as
natural transformations and drawing the internal views of their
internal views:
%
$$\pu
  \diag{adj-nat-conditions}
$$

%
\subsection{A way to teach adjunctions \DONE}
\label{teaching-adjunctions}

I mentioned in the first sections that I have tested some parts of
this language in a seminar course --- described here:
\cite{OchsWLD2019} --- and that in it I teach Categories starting by
adjunctions. Here's how: we start by the basics of $λ$-calculus and
some sections of \cite{PH1}, and then I ask the students to define
each one of the operations in the right half of the diagram below as
$λ$-terms:
$$
  \pu
  \diag{generic-adjunction-big}
  \qquad
  \diag{(xB)-|(B->)}
$$

Then we see the definition of functors, natural transformations and
adjunctions, and we check that the right half is a particular case of
the diagram for a generic adjunction in the left half. After that, and
after also checking that in the Planar Heyting Algebras of \cite{PH1}
we have an adjunction $(∧Q)⊣(Q→)$, I help the students to decypher
some excerpts of standard texts on CT --- in the last time that I gave
the course we used \cite{Awodey}, but I am planning to use \cite{CWM2}
the next time.

\msk

From the components of the generic adjunction in the diagram above it
is possible to build this big diagram:
$$\pu
  \diag{teaching-adjunctions-1}
$$

Let's use these names for its subdiagrams: $\sm{ A \\ BCDEF \\ G \\ I}$.

A {\sl fully-specified adjunction} between categories $\catB$ and
$\catA$ has lots of components: $(L, R, ε, η, ♭, ♯, \univ(ε),
\univ(η))$, and maybe even others, but usually we define only some of
these components; there is a Big Theorem About Adjunctions (below!)
that says how to reconstruct the fully-specified adjunction from some
of its components.

Some parts of the diagram above can be interpreted as definitions,
like these:
$$\begin{array}{c}
  Lf := (η_A∘f)^♭ \\
  [5pt]
  g := ε_B∘Lh 
    \qquad ε_B := (\id_{RB})^♭
    \qquad η_A := (\id_{LA})^♯
    \qquad h := Rg∘η_A \\
  [5pt]
  Rk := (k∘η_B)^♯ \\
  \end{array}
$$

The subdiagrams $B$ and $F$ can also be interpreted in the opposite
direction, as:
$$\begin{array}{rclcrcl}
  g^♯ &:=& (∀A.∀g.∃!h)Ag    &\phantom{mmmmmm}&  h^♭ &:=& (∀B.∀h.∃!g)Bh \\
       &=& (\univ_{ε_B})Ag  &&                       &=& (\univ_{η_A}) Bh \\
  \end{array}
$$

The notations $(∀A.∀g.∃!h)Ag$ and $(\univ_{ε_B})Ag$ are clearly abuses
of language --- but it's not hard to translate them to something
formal, and they inspire great discussions in the classroom... also,
they can help us to understand and formalize constructions like this
one,
%
$$\pu
  Lf := (\univ_{η_A})(LA)(η_A∘f)
  \qquad
  \diag{building-L_1f}
$$
that are needed in cases like the part (ii) of the Big Theorem.

The Big Theorem About Adjunctions is this --- it's the Theorem 2 in
\cite[page 83]{CWM2}, but with letters changed to match the ones we
are using in our diagrams:

\def\ORIG#1{\msk\ColorBrown{#1}}
\def\ORIG#1{}

\newpage

\begin{quotation}

  \ORIG{{\bf Theorem 2.} Each adjunction $〈F,G,φ〉: X \rightharpoonup
    A$ is completely determined by the items in any one of the
    following lists:}

  {\bf Big Theorem About Adjunctions.} Each adjunction $〈L,R,♯〉: \catA
  \rightharpoonup \catB$ is completely determined by the items in any
  one of the following lists:

  \ORIG{(i) Functors $F$, $G$, and a natural transformation $η: 1_X
    \tnto GF$ such that each $η_x: x→GFx$ is universal to $G$ from
    $x$. Then $φ$ is defined by (6).}

  (i) Functors $L$, $R$, and a natural transformation $η:
  \id_\catA→RL$ such that each $η_A: A→RLA$ is universal to $R$ from
  $A$. Then $♯$ is defined by (6).

  \ORIG{(ii) The functor $G: A → X$ and for each $x∈X$ an
    object $F_0x∈A$ and a universal arrow $η_x:x→GF_0x$ from $x$
    to $G$. Then the functor $F$ has object function $F_0$ and is
    defined on arrows $h:x→x'$ by $GFh∘η_x = η_{x'}∘h$.}

  (ii) The functor $R: \catB → \catA$ and for each $A∈\catA$ an object
  $L_0A∈\catB$ and a universal arrow $η_A:A→RL_0A$ from $A$ to $R$.
  Then the functor $L$ has object function $L_0$ and is defined on
  arrows $f:A'→A$ by $RLf∘η_{A'} = η_A∘f$.

  \ORIG{(iii) Functors $F$, $G$, and a natural transformation $ε: FG
    \tnto I_A$ such that each $ε_a:FGa→a$ is universal from $F$ to
    $a$. Here $φ^{-1}$ is defined by (7).}

  (iii) Functors $L$, $R$, and a natural transformation $ε:
  LR→\id_\catB$ such that each $ε_B:LRB→B$ is universal from $L$ to
  $B$. Here $♭$ is defined by (7).

  \ORIG{(iv) The functor $F:X→A$ and for each $a∈A$ an object $G_0a∈X$
    and an arrow $ε_a:FG_0a→a$ universal from $F$ to $a$.}

  (iv) The functor $L:\catA→\catB$ and for each $B∈\catB$ an object
  $R_0B∈\catA$ and an arrow $ε_B:LR_0B→B$ universal from $L$ to $B$.

  \ORIG{(v) Functors $F$, $G$ and natural transformations $η:I_x \tnto
    GF$ and $ε: FG \tnto I_A$ such that both composites (8) are the
    identity transformations. Here $φ$ is defined by (6) and $φ^{-1}$
    by (7).}

  (v) Functors $L$, $R$ and natural transformations $η:\id_\catA→RL$
  and $ε:LR→\id_\catB$ such that both composites (8) are the identity
  transformations. Here $♯$ is defined by (6) and $♭$ by (7).

\end{quotation}

My plan for the next incarnation of the course is to ask the students
to 1) visualize in the big diagram all the objects and constructions
in the Big Theorem, 2) take the original Theorem 2 in \cite{CWM2} and
draw the missing diagrams for it, 3) decypher some other parts of the
section about adjunctions in \cite{CWM2}.


\newpage

%

\section{The Basic Example as a skeleton \DONE}
\label{basic-example-as-skel}

In the sections \ref{the-conventions} and \ref{to-deserve-a-name} I
claimed that the diagram of the Basic Example is a ``skeleton'' of a
certain theorem, in the sense that both the statement and the proof of
that theorem can be reconstructed from just the diagram and very few
extra hints. Let's see the details of this.

%
\subsection{Reconstructing its functors \DONE}
\label{basic-example-functors}

\def\Yzero    {\mathsf{Y0}}
\def\Yzeroplus{\mathsf{Y0^+}}
\def\Yone     {\mathsf{Y1}}

Let's call this diagram --- the diagram of the Basic Example ---
$\Yzero$:
$$\Yzero \qquad := \quad \diag{Basic-Example}$$

We don't know yet the precise meaning of the functors $\catB(C,-)$ and
$\catA(A,R-)$, but if we enlarge $\Yzero$ to
%
%
$$\pu
  \Yzeroplus \qquad := \quad \diag{Basic-Example-plus}
$$
and we draw the internal views of $\catB(C,-)$ and $\catA(A,R-)$ then
the meanings for $\catB(C,-)$ and $\catA(A,R-)$ become obvious:
%
$$\pu
  \diag{basic-example-obvious-functors}
$$

So:
$$\begin{array}{rcl}
  \catB(C,-)   &:&  \catB → \Set \\
  \catB(C,-)_0 &:=& λD.\catB(C,D) \\
  \catB(C,-)_1 &:=& λg.λf.g∘f \\
  [5pt]
  \catA(A,R-)   &:&  \catB → \Set \\
  \catA(A,R-)_0 &:=& λD.\catA(A,RD) \\
  \catA(A,R-)_1 &:=& λg.λh.Rg∘h \\
  \end{array}
$$

%
\subsection{Reconstructing its natural transformation \DONE}
\label{basic-example-NT}

We also don't know --- yet --- what is the natural transformation
$$\catB(C,-) \ton{T} \catA(A,R-).$$
Its internal view is this:
%
$$\pu
  \diag{basic-example-obvious-NT}
$$
Note that we only drew the vertical arrows of the internal view of the internal view.

If we have an arrow $η:A→RC$ then we have a natural construction for
$T_0$: $TD(f):=Rf∘η$, and we can redraw the internal view of the
internal view as:
%
$$\pu
  \diag{basic-example-obvious-NT-2}
$$
The square condition clearly holds, because:
$$\begin{array}{rcl}
  Rg∘(Rf∘η) &=& (Rg∘Rf)∘η \\
            &=& R(g∘f)∘η. \\
  \end{array}
$$

So
$$\begin{array}{rcl}
    T_0 &:=& λD.λf.Rf∘η. \\
  \end{array}
$$

%
\subsection{Reconstructing its bijection \DONE}
\label{basic-example-bij}

We can give names like `$d$' and `$u$' for the two components of the
curved bijection, like this:
%
$$\pu
  \diag{basic-example-bij-1}
$$
but the notation at the right will be clearer.

We just saw how the direction `$d$' of the bijection works:
$$\begin{array}{rcl}
    (T_η)_0 &:=& λD.λf.Rf∘η. \\
  \end{array}
$$

Here's how to find a natural construction for $u$. Suppose that we
have a natural transformation $T$. Then $TC(\id_C)$ is an element of
$\catA(A,RC)$:
%
$$\pu
  \diag{basic-example-bij-2}
$$

We can define:
$$\begin{array}{rcl}
    η_T &:=& TC(\id_C). \\
  \end{array}
$$

Now we need to check that $d$ and $u$ are mutually inverse, or, in the
other notation, that the round trips $η \mapsto T_η \mapsto η_{(T_η)}$
and $T \mapsto η_T \mapsto T_{(η_T)}$ are identity maps. Here is a
good way to draw the round trips:
%
$$\pu
  \diag{basic-example-bij-3}
$$

Checking that $η \mapsto T_η \mapsto η_{(T_η)}$ yields back the
original $η$ is easy --- we just have to start with $η_{(T_η)}$ and
reduce it as most as we can:
%
%
$$\begin{array}{rcl}
 η_{(T_η)} &=& T_ηC(\id_C) \\
           &=& (λD.λg.(Rg∘η)) C(\id_C) \\
           &=& (λg.(Rg∘η)) (\id_C) \\
           &=& R(\id_C)∘η \\
           &=& \id_{RC}∘η \\
           &=& η \\
 \end{array}
$$

Checking that the other round trip, $T \mapsto η_T \mapsto T_{(η_T)}$,
yields back the original $T$ is not trivial. In the terminology of the
convention (CSk) from Section \ref{the-conventions}, to reconstruct
that proof we need an extra hint: that at some point in the proof we
will have to use that the original $T$ obeys $\sqcond_T$, and that
we will have to ``evaluate'' $\sqcond_T$ on these inputs:

$$\pu
  \diag{Y0-NT-2}
$$

This yields:
%
$$\pu
  \diag{Y0-NT-3}
$$
so $Rf∘(TC(\id_C)) = TD(f)$.

We want to check that for all $D$ and $f$ we have $T_{(η_T)}D(f) =
TD(f)$. We have:
$$\begin{array}{rcl}
  T_{(η_T)}D(f) &=& (λD.λf.Rf∘η_T)D(f) \\
                &=& (λf.Rf∘η_T)(f) \\
                &=& Rf∘η_T \\
                &=& Rf∘(TC(\id_C)) \\
                &=& TD(f). \\
  \end{array}
$$

It works! So we have a natural construction for the bijection $T ↔ η$,
given by:
$$\begin{array}{rcl}
  T_0 &:=& λD.λf.Rf∘η \\
    η &:=& TC(\id_C) \\
  \end{array}
$$

%
\subsection{The full reconstruction \DONE}
\label{basic-example-full}

We have just reconstructed all the typings and definitions for the
diagram $\Yzero$. Here is the full reconstruction, except for the
``proof terms'' like $\respids$, $\assoc$, $\idL$ and $\idR$ for each
functor, $\sqcond$ for each natural transformations, and the proofs
that both round trips in the bijections are identity maps:
$$\diag{Basic-Example}
  \qquad
  \begin{array}{rl}
    & \catA \text{ is a category}, \\
    & \catB \text{ is a category}, \\
    & R:\catB \to \catA, \\
    & A ∈ \catA, \\
    & C ∈ \catB, \\
    & η:A→RD, \\
      [5pt]
    & \catB(C,-)    :  \catB → \Set,   \\
    & \catB(C,-)_0  := λD.\catB(C,D),  \\
    & \catB(C,-)_1  := λg.λf.g∘f,      \\
      [5pt]
    & \catA(A,R-)   : \catA → \Set,    \\
    & \catA(A,R-)_0 := λD.\catA(A,RD), \\
    & \catA(A,R-)_1 := λg.λh.Rg∘h,     \\
      [5pt]
    & T : \catB(C,-) → \catA(A,R-), \\
      [5pt]
    & T_0 := λD.λf.Rf∘η, \\
    & η := TC(\id_C). \\
  \end{array}
$$


It shouldn't be hard --- for someone with practice --- to translate
the types and definitions at the right above to the language of some
proof assistant. I tried to do this in Idris (\cite{Brady}) using
\cite{IdrisCT2019} but I didn't go very far... I implemented the
protocategories, protofunctors and proto-NTs of \cite[section
  19]{IDARCT} to be able to skip the proof terms on my first
prototypes, but I got stuck trying to implement the formalization of
$\Yzero$ as a single datatype...

\bsk

{\sl (Help would be greatly appreciated!...)}

\newpage

%
\section{Extensions to the diagrammatic language \DONE}
\label{extensions}

Our diagrammatic language and the list of conventions in Section
\ref{the-conventions} can be extended --- ``by the user'' --- in
zillions of ways. Let's see some examples of extensions.

%
\subsection{A way to define new categories \DONE}
\label{comma-categories}

We saw in the sections \ref{internal-view-functor} and
\ref{basic-example-functors} how to use diagrams to define functors,
and in sections \ref{internal-view-NT} and \ref{basic-example-NT} how
to define natural transformations. We can define new categories by
diagrams, too.

\pu
\def\commaobj#1#2#3#4{{
  \left(        \def\A{#1}
                \def\f{#4}
    \def\B{#2} \def\FB{#3}
     \diag{comma-obj-0}
  \right)
  }}


\def\dnAR{{(A{↓}R)}}

$$\pu
  \diag{defining-a-comma-category}
$$

\def\AProofOf   #1{\llangle#1\rrangle}
\def\AllProofsOf#1{\llbracket#1\rrbracket}

My favorite way --- a syntax sugar! --- of visualizing the comma
category $\dnAR$ is the middle third of the diagram above, in which
the objects of $\dnAR$ are depicted as L-shaped diagrams. To
understand the typings and the commutativity conditions we have to
look at the left third --- it indicates that $f$ must obey $Rf∘η=g$.
The right third shows a generic morphism in $\dnAR$ without the syntax
sugar, but we still have to look at the left third to type it. We
have:
$$\begin{array}{rl}
  \text{In a context in which}
    & \catA \text{ is a category}, \\
    & \catB \text{ is a category}, \\
    & R : \catB → \catA, \\
    & A \text{ is an object of $\catA$}, \\
  \text{we define the category}
    & \dnAR \text{ as follows:} \\
  [5pt]
  \text{An object of}
    & \dnAR \\
  \text{is a pair}
    & (C,η) \\
  \text{in which}
    & C : \catB_0 \\
  \text{and}
    & η : \Hom_\catA(A,RC); \\
  \text{so}
    & (C,η) : ΣC⠆\catB_0. \Hom_\catA(A,RC) \\
  \text{and}
    & \dnAR_0 := ΣC⠆\catB_0. \Hom_\catA(A,RC). \\
  [5pt]
  \text{A morphism}
    & f: (C,η) → (D,g) \text{ in $\dnAR$} \\
  \text{is an}
    & f: \Hom_\catB(C,D) \text{ such that $Rf∘η=g$}, \\
  \text{or equivalently a pair}
    & (f,\AProofOf{Rf∘η=g}); \\
  \text{we have}
    & (f,\AProofOf{Rf∘η=g}) : Σf⠆\Hom_\catB(C,D).\AllProofsOf{Rf∘η=g}, \\
  \text{so}
    & \Hom_\dnAR((C,η),(D,g)) := \\
    & Σf⠆\Hom_\catB(C,D).\AllProofsOf{Rf∘η=g}.
  \end{array}
$$

The notations $\AProofOf{P}$ and $\AllProofsOf{P}$ are non-standard.
For any proposition $P$ we denote by $\AllProofsOf{P}$ the set of
witnesses of $P$ (see \cite[p.18]{HOTT}) and by $\AProofOf{P}$ a
witness that $P$ is true; formally, $\AProofOf{P}$ is a variable (with
a long name!) whose type is $\AllProofsOf{P}$, and $\AllProofsOf{P}$
is a singleton when $P$ is true and the empty set when $P$ is false. A
good way to remember this notation is that $\AllProofsOf{P}$ looks
like a box and $\AProofOf{P}$ looks like something that comes in that
box.

\msk

This defines formally the first two components of the category
$\dnAR$. Remember that a category $\catC$ has seven components:
$$\catC = (\catC_0, \Hom_\catC, \id_\catC, ∘_\catC;
   \assoc_\catC, \idL_\catC, \idR_\catC)
$$
We are pretending that the other components of $\dnAR$ are ``obvious''
in the sense of Section \ref{to-deserve-a-name}.


%
\subsection{The Yoneda Lemma \DONE}
\label{yoneda-lemma}

The formalization of $\Yzero$ as a series of typings and definitions
in Section \ref{basic-example-full} {\sl suggests} that {\sl some}
operations from Type Theory that can be applied on the formalization
side should be translatable to the diagram side; for example,
substitution. This one clearly works: if we substitute $\catA$ by
$\Set$ and $A$ by the set 1 we get this,
%
%
\pu
$$
  \Yzero \bmat{\catA := \Set \\ A := 1}
  \qquad = \quad
  \diag{Basic-Example-using-Set-and-1}
$$

For each $D∈\catB$ we have a bijection $\Set(1,RD) ↔ RD$ --- and we
can use these bijections to build a natural isomorphism $\Set(1,R-) ↔
R$, that we will add to the diagram:

%
\pu
$$
  \Yone
  \qquad := \qquad
  \diag{Basic-Example-using-Set-and-1-and-R}
$$

We can obtain $T'$ from $T$ and vice-versa by composing them with
$\Set(1,R-) ↔ R$.

The diagram $\Yone$ ``is'' the Yoneda Lemma --- but it doesn't have a
single top-level arrow, so we can't apply the convention (CTL) to it,
and we need to specify its ``meaning'' explicitly. The statement of
the Yoneda Lemma is that there is a bijection
$$RC ↔ \Hom(\catB(C,-),R);$$
Once we know that it is easy to see that the diagram $\Yone$ shows how
we can build it by combining three bijections that we understand well:
$$\begin{array}{l}
  RC \\
  ↔ \Hom(1,RC) \\
  ↔ \Hom(\catB(C,-),\Set(1,R-)) \\
  ↔ \Hom(\catB(C,-),R) \\
  \end{array}
$$
So $\Yone$ shows a way to build the bijection $RC ↔
\Hom(\catB(C,-),R)$.

%
\subsection{The Yoneda embedding \DONE}
\label{yoneda-embedding}

Let $B$ be an object of $\catB$. If we replace the functor $R:
\catB→\Set$ in $\Yone$ by $\catB(B,-)$ and do some other renamings we
get this:
%
%
\pu
$$
  \Yone\bmat{
    R := \catB(B,-) \\
    η := \nameof{f} \\
    T := T' \\
    T' := T \\
    }
  \quad := \quad
  \diag{Basic-Example-using-Hom(B,-)}
$$

We can consider that the diagram above is a skeleton for the {\sl
  proof} that there is a bijection between arrows $f:B→C$ and natural
transformations $T:\catB(C,-)→\catB(B,-)$. The two directions of the
bijection are easy to define, as $T_0 := λD.λg.g∘f$ and $f :=
TC(\id_C)$, but the proof that the round trips $f \mapsto T \mapsto f$
and $T \mapsto f \mapsto T$ give back the original $f$ and $T$ are
tricky, as we saw in Section \ref{basic-example-bij}.

Usually people draw a simple diagram that just {\sl states} that the
obvious map $\catB(B,C) → \Hom(\catB(C,-),\catB(B,-)$ is a bijection,
somehow like this:
%
%
$$\pu
  \diag{Y-emb-square-only}
$$
Compare with \cite[p.60]{Riehl}; note that our arrow in the middle of
the square is a `$↔$'.

We can draw it with more details as:
%
$$\pu
  \diag{Y-emb-full}
$$
Note that it defines a {\sl contravariant} functor $𝐛y:
\catB^\op→\Set$ whose action on objects is $C \mapsto \catB(C,-)$.

We consider that the morphism $f:B→C$ in the diagram is inside
$\catB$, not inside $\catB^\op$. This is explained in the next
section.

%
\subsection{Opposite categories \DONE}
\label{opposite-categories}



\def\Aop{{\catA^\op}}

Suppose that we have a diagram $A \ton{f} B \ton{g} C$ in a category
$\catA$. There are several different notations for the corresponding
diagram in $\Aop$: for example, in \cite[p.33]{CWM2} it would be
written as $A \otn{f^\op} B \otn{g^\op} C$, while in
\cite[p.15]{AbramskyTzevelekos} as $A \otn{f} B \otn{g} C$. The
convention (COT) says that the notation in our diagrams should be as
close as possible to the notation in the original text --- so let's
see how to support the notation in \cite{AbramskyTzevelekos}, that
looks a bit harder than the one in \cite{CWM2}.

We want to define a new category, $\Aop$, using tricks similar to the
ones in Section \ref{comma-categories}, but now we can't pretend that
the new composition is obvious. We will define $(\Aop)_0$,
$\Hom_\Aop$, $\id_\Aop$, and $∘_\Aop$ without any textual
explanations, with just the diagrams to convince the readed that our
definitions are reasonable.
%
$$\pu
  \diag{A-and-Aop}
  \quad
  \begin{array}{c}
  \catA_0 =: (\Aop)_0 \\
  \\
  \Hom_\catA(A,B) =: \Hom_\Aop(B,A) \\
  \\
  \id_\catA(A) =: \id_\Aop(A) \\
  \\
  g ∘_\catA f =: f ∘_\Aop g \\
  \\
  \\
  \\
  \\
  \\
  \end{array}
$$

In the diagram below $F:\catA^\op→\catB$ is a contravariant functor,
and the $\catA$ above $\catA^\op$ indicates that $g:C→D$ is a morphism
of $\catA$, not of $\catA^\op$. I am not very happy with this trick
but I haven't found a better alternative yet.
%
$$\pu
  \diag{contravariant-functor}
$$

\msk

%
\subsection{Universalness as something extra \DONE}
\label{ness}

We can consider that an universal arrow is an arrow $η:A→RC$ with an
extra property; I showed at the end of Section
\ref{freyd-with-functors} how to think of that property as being just
$∀D.∀g.∃!f$, and how to treat that as an abbreviation for something
bigger and more formal.

We can also treat a universal arrow as an arrow $η:A→RC$ plus extra
{\sl structure} --- this extra structure is an operation that returns
for each $D$ an inverse for the operation $g \mapsto Rg∘η$. For more
on properties and structure, see \cite[p.15]{BaezShulman2007}.

In any case this ``universalness'' can be treated as
something extra, and a universal arrow can be expressed as:
$$(η,\univ_η)$$
using dependent types.

Several of these ``-ness''es have standard graphical representations:
for example pullbackness is indicated by a `$\pbsymbol7$', and
monicness is indicated by a tail like this: `$\monicto$'.
\cite{FreydScedrov} defines lots of graphical representations for
``-ness''es starting on its page 37. We will use an `$:=$' to define a
new annotation that is an abbreviation for extra structure:
%
$$\pu
  \diag{universalness}
$$

This is pullbackness:
%
$$\pu
  \diag{pullbackness}
$$

%
\subsection{Representable functors \DONE}
\label{representable-functors}

It is easy to see that in $\Yzero$ the universality of $η$ is
equivalent to the natural-iso-ness of $T$; in $\Yone$ the universality
of $η$ is equivalent to the natural-iso-ness of $T$, and this is
equivalent to the natural-iso-ness of $T'$. The constructions should
be evident from these diagrams:
%
%
$$\pu
  \diag{univ-arrows-univ-elts}
$$

The diagram at the right above can be seen as the missing diagram for
Proposition 2 in \cite[p.60]{CWM2}, that says this (I've translated
its letters to the ones I use):


\begin{quotation}

  {\bf Definition.} {\sl Let $\catB$ have small hom-sets. A
    representation of a functor $R:\catB→\Set$ is a pair $〈C,T'〉$,
    with $C$ an object of $\catB$ and
    $$T':\catB(C,-)→R$$
    a natural isomorphism. The object $C$ is called the representing
    object. The functor $R$ is said to be representable when such a
    representation exists.}

  \msk

  Up to isomorphism, a representable functor is thus just a covariant
  hom-functor $\catB(C,-)$. This notion can be related to universal arrows as
  follows.

  \msk

  {\bf Proposition 2.} Let 1 denote any one-point set and let $\catB$
  have small hom-sets. If $〈C, η:1→RC〉$ is a universal arrow from $1$
  to $R: \catB→\Set$, then the function $T'$ which for each object $D$
  of $\catB$ sends the arrow $f:C→D$ to $(Rf)(η(*))∈RD$ is a
  representation of $R$. Every representation of $R$ is obtained in
  this way from exactly one such universal arrow.

\end{quotation}

The operations $T'↦η$ and $η↦T'$ can be defined as:
$$\begin{array}{rcl}
  η  &:&  1→RC             \\
  T' &:&  \catB(C,-)→R     \\
  η  &:=& λ*.(T'C(\id(C))) \\
  T' &:=& λD.λf.(Rf)(η(*)) \\
  \end{array}
$$

%
\subsection{An example of a representable functor \DONE}
\label{representable-functor-ex}


Emily Riehl gives two pages of examples of representable functors in
\cite[pages 51--53]{Riehl}. Her example (iv) is:

\begin{quotation}
\begin{enumerate}

\item[(iv)] The functor $U:\Ring→\Set$ is represented by the unital
  ring $\Z[x]$, the polynomial ring in one variable with integer
  coefficients. A unital ring homomorphism $Z[x]→R$ is uniquely
  determined by the image of $x$; put another way, $\Z[x]$ is the {\sl
    free unital ring on a single generator}.

\end{enumerate}
\end{quotation}


She develops more this example in page 63, as:

\begin{quotation}

  {\bf Example 2.3.4.} Recall from Example 2.1.5(iv) that the
  forgetful functor $U:\Ring→\Set$ is represented by the ring $Z[x]$.
  The universal element, which defines the natural isomorphism
  $$\Ring(Z[x], R) ≅ UR,$$
  is the element $x∈\Z[x]$. As in the proof of the Yoneda lemma, the
  bijection above is implemented by evaluating a ring homomorphism
  $ϕ:\Z[x]→R$ at the element $x∈\Z[x]$ to obtain an element $ϕ(x)∈R$.

\end{quotation}

Here is the ``missing diagram'' for both excerpts:
%
$$\pu
  \diag{2.3.4._Z[x]}
$$

That diagram may be a good starting point to explain the Yoneda Lemma
to ``children''.


%

\subsection{Functors as objects \DONE}
\label{functors-as-objects}

One way to treat a diagram in $\Set$ like this
%
\pu
$$
  F \qquad := \qquad \diag{evil-presheaf}
$$
as a functor is to think that that diagram is an abbreviation --- it
is just the upper-right part of a diagram like this,
%
\pu
%
$$\pu
  \diag{evil-presheaf-as-functor}
$$
where we add the extra hint that the index category $\catK$ is exactly
the kite-shaped preorder category drawn above the ``$\catK$''.

The convention (CFSh) says that the image by a functor of a diagram is
a diagram with the same shape, so according to that convention we have
$F(1) = \{24,25\}$, $F(4→5) = (\{1\} \ton{1↦1} \{0,1\})$, and so
on; so the upper right part of the diagram above {\sl defines} $F$.

Note that the single `$↦$' above the $\catK \ton{F} \Set$ stands for
several `$↦$'s, one for each object and one for each morphism, and
note that $F$ is an object of $\Set^\catK$.

%
\subsection{Geometric morphisms for children \DONE}
\label{gms-for-children}

Let $\catA$ and $\catB$ be these preorder categories, and let
$f:\catA→\catB$ be the inclusion functor from $\catA$ to $\catB$:
$$
  A := \pshAargs2345
  \qquad
  B := \pshBargs123456
$$

The left half of the diagram below is the standard definition of a
geometric morphism $f$ from a topos $\calE$ to a topos $\calF$. A
geometric morphism $f:\calE→\calF$ is actually an adjunction $f^*⊣f_*$
plus the guarantee that $f^*:\calE \ot \calF$ preserves limits, which
is a condition slightly weaker than requiring that $f^*$ has a left
adjoint. When that left adjoint exists it is denoted by $f^!$, and we
say that $f^!⊣f^*⊣f_*$ is an {\sl essential geometric morphism}. The
only non-standard thing about the diagram at the left below is that is
contains an internal view of the adjunction $f^*⊣f_*$.
$$
  \diag{gm-for-adults}
  \qquad
  \def\LG{\pshAargs{G_2}{G_3}{G_4}{G_5}}
  \def\G {\pshBargs{G_1}{G_2}{G_3}{G_4}{G_5}{G_6}}
  \def\H {\pshAargs{H_2}{H_3}{H_4}{H_5}}
  \def\RH{\pshBargs{H_2{×_{H_4}}H_3}{H_2}{H_3}{H_4}{H_5}{1}}
  \diag{gm-for-children}
$$

The right half of the diagram is a particular case of the left half.
Its lower line, $\catA \ton{f} \catB$, does not exist in the left
half. The inclusion functor $f$ induces adjunctions $f^!⊣f^*⊣f_*$ as
this,
%
$$\pu
  \diag{essential-GM-small}
$$
where $f^*$ is easy to define and $f^!$ and $f_*$ not so much --- the
standard way to define $f^!$ and $f_*$ is by Kan extensions.

The big square in the upper part of the diagram is an internal view of
the adjunction $f^*⊣f_*$, with the functors $f^*G$, $G$, $H$, and
$f_*H$ being displayed as their internal views. We can choose the sets
$G_1, \ldots, G_6$ and the morphisms between them arbitrarily, so this
is an internal view of an arbitrary functor $G:\catB→\Set$; and the
same for $H$.

The arrow $f^*G \mapsot G$ can be read as a definition for the action
of $f^*$ on objects --- it just erases some parts of the diagram ---
and the arrow $H \mapsto f_*H$ can be read as a definition for the
action of $f_*$ on objects --- $f_*$ ``reconstructs'' $H_1$ and $H_6$
in a certain natural way. It is easy to reconstruct the actions of
$f^*$ and $f_*$ on morphisms from just what is shown, and to
reconstruct the two directions of the bijection.

The big diagram above can be used 1) to convince people that are not
hardcore toposophers that this diagrammatic language can make some
difficult categorical concepts more accessible, and 2) as a starting
point to generate diagrams ``for children'' for several parts of the
Elephant, and even to prove new theorems on toposes. For more on (1),
see \cite{OchsLucatelli} and \cite{OchsVGMS2018}; for (2), see
\cite{MDE}.

%
\subsection{Reading the Elephant \DONE}
\label{reading-the-elephant}

In Section \ref{teaching-adjunctions} we saw a strategy for helping
(beginner) students to read a difficult text on CT: we start with
diagrams for the most important concepts, in both a general case ``for
adults'' and a well-chosen particular case ``for children'', we give
them exercises to make sure that they understand the constructions in
the case ``for children'', we give them a few more exercises to make
sure that they understand the general case, we ask them to read
excerpts from a standard textbook in a version where the letters were
changed to match the diagrams, and then we ask them to work on the
original version of these excerpts with the original notation, and on
some other parts of the same chapter... this can be done for the
Elephant too --- here are the parts that are more relevant for our
diagrams on geometric morphisms, with the notation adjusted:

\begin{quotation}


  {\bf Definition 4.1.1.} (a) Let $\calF$ and $\calE$ be toposes. A
  geometric morphism $f: \calE → \calF$ consists of a pair of functors
  $f_*: \calE → \calF$ (the direct image of f) and $f^*: \calF →
  \calE$ (the inverse image of $f$) together with an adjunction ($f^*
  ⊣ f_*$), such that $f^*$ is cartesian (i.e. preserves finite
  limits).

  \msk

  (...)

  \msk


  {\bf Example 4.1.4.} Let $f: \catA → \catB$ be a functor between
  small categories. Then composition with $f$ defines a functor $f^*:
  \Set^\catB→\Set^\catA$, which has adjoints on both sides, the left
  and right {\sl Kan extensions} along $f$: for example, the right Kan
  extension $\Ran_f$ sends a functor $H: \Set^\catA$ to the functor
  whose value at an object $B$ of $\catB$ is the limit of the diagram
  %
  $$ (B ↓ f) \diagxyto/->/^U \catA \diagxyto/->/^H \Set $$
  (here $(B ↓ f)$ is the comma category whose objects are pairs
  $(A,ϕ)$ with $ϕ:B→fA$ in $\catB$, and $U$ is the forgetful functor
  from this category to $\catA$). Thus $f^*$ is the inverse image of a
  geometric morphism $\Set^\catA → \Set^\catB$, whose direct image is
  $\Ran_f$.

  \msk

  (...)

  \msk


  We note that the geometric morphisms which arise as in 4.1.4, though
  not as special as those of 4.1.2, still have the property that their
  inverse image functors have left adjoints as well as right adjoints.
  We call a geometric morphism $f$ {\it essential} if it has this
  property; we normally write $f_!$ for the left adjoint of $f^*$.
  With the aid of this notion, we can prove a partial converse to
  4.1.4:

  \ssk


  {\bf Lemma 4.1.5.} Let $\catA$ and $\catB$ be small categories such
  that $\calB$ is Cauchy-complete (cf.\ 1.1.10). Then every essential
  geometric morphism $f: \Set^\catA → \Set^\catB$ is induced as in
  4.1.4 by a functor $\catA → \catB$.

  \msk

  (...)

  \msk


  {\bf Proposition 4.2.8.} With the notation established above, the counit
  $h^*h_*→1$ is an isomorphism.

  (...)

  A geometric morphism $h$ satisfying the condition that the counit
  $h^*h_*→1$ is an isomorphism, or the equivalent condition that $h_*$
  is full and faithful, is called an {\sl inclusion} (though some
  authors prefer the term {\sl embedding}). We shall study inclusions
  in greater detail in the next three sections; for the present, we
  digress briefly to note an alternative characterization of them:


\end{quotation}

The really interesting part would be to show that the unit $η$ of the
adjunction $f^*⊣f_*$ ``is'' a sheafification functor, and that the
geometric morphism for children of the diagram yields an example of
sheaf... but that would need lots of different fragments from several
different sections of the book.

%
\section{How to name this diagrammatic language \DONE}
\label{how-to-name-this-language}

{\sl I don't have any idea!...}

It can be used to produce missing diagrams, and sometimes these
missing diagrams are skeletons. We can use it to work in two styles in
parallel, ``for adults'' and ``for children''... maybe something like
``Missing Skeletons for Children''?

\msk

Suggestions welcome.

%

\section{Why ``my conventions''? \DONE}
\label{why-my-conventions}

I learned CT as an autodidact in a totally disorganized way. In the
first years I just read, or rather {\sl tried} to read, everything
that was available in my university's library, trying to locate the
parts that could be useful to my main interest at that time, that was
Non-Standard Analysis and how to do something similar to NSA but using
filter-powers instead of ultrapowers...

It was only after that that I realized that I had to learn how to {\sl
  write}. I remember one time spending a whole evening on an exercise
of the beginning of \cite{LambekScott} that says just ``prove that for
categories $\catA$, $\catB$, and $\catC$ we have $\catA^{\catB×\catC}
≅ (\catA^\catB)^\catC$'' --- the full proof had lots of parts, and I
saw that I didn't know how to organize them in a neat way... also, the
proofs given in books and articles just state the main parts and
pretend that the rest is obvious, and in the case of
$\catA^{\catB×\catC} ≅ (\catA^\catB)^\catC$ there were no ``main
parts'', so I had to learn how to write down a proof in full, and this
was a new style to me...

Even now, many years after that, I still have the sensation that I had
to improvise practically everything in my ways --- both the
``algebraic'' way and the ``diagrammatic'' way --- of writing
categorical proofs, and that I still don't know even a tiny fraction
of the techniques for writing that people learn when they take CT
courses and they have opportunities to discuss exercises with other
students and with TAs and more senior people...

\msk

{\sl The ``my conventions'' in the title of this text, and my use of
  the first person everywhere, are a way to stress that I still don't
  know enough about other people's private languages for CT, and that
  this is an attempt to gain access to other private languages,
  diagrammatic or not... I am especially interested in how people
  write when they turn their level-of-detail knob to a very high
  position.}


%
\section{Related and unrelated work \DONE}
\label{related-and-unrelated}


The diagrammatic language that I described here seems to be unrelated
to the ones in \cite{CoeckePQP} and \cite{CoeckeNewStruP} --- that
describe {\sl lots} of diagrammatic languages --- and also unrelated to
\cite{MarsdenCTUSD}.


I've taken an approach that is the opposite of \cite{CaccamoWinskel}.
Cáccamo and Winskel define a derivation system that can only construct
functors, natural transformations, etc, that obey the expected
naturality conditions, while we allow some kinds of sloppinesses, like
constructing something that looks like a functor and pretending that
it is a functor when it may not be. When I started working on this
diagrammatic language I had a companion derivation system for it;
\cite[Section 14]{IDARCT} mentions it briefly, but it doesn't show the
introduction rules that create (proto)functors and (proto)natural
transformations and that allow being sloppy (``in the syntactical
world'').


Some of my excuses for allowing one to pretend that a functor is a
functor and leaving the verification to a second stage come from
\cite{ChengMorally}. I learned a {\sl lot} on how mathematicians use
intuition and diagrams from \cite{Kromer} --- \cite{KromerSlides} is a
great summary --- and \cite{Corfield}, and they have helped me to
identify which characteristics of my diagrammatic language are very
unusual and may be new, and that deserve to be presented in detail.


Many of the first ideas for my diagrammatic language appeared when I
was reading \cite{SeelyBeck}, \cite{SeelyLCCC}, \cite{SeelyPLC},
\cite{Jacobs}, and \cite{SeelyDiff} and trying to draw the ``missing
diagrams'' in those papers in both the original notation and in the
``archetypal case'' (\cite[Section 16]{IDARCT}).


%

\section{What next? \DONE}
\label{what-next}

At this point I think that it is more interesting to ``implement''
more categorical definitions and proofs in this diagrammatic language
than to try to formalize it completely or try to prove meta-theorems
about it. I am doing that by (re)reading parts of several papers and
articles and drawing the missing diagrams in them; for details and
links, see:

\msk

\centerline{\url{http://angg.twu.net/math-b.html\#favorite-conventions}}

\msk

Besides this, here's what I've planned for the next steps. Most of
them can be done in parallel.

\begin{enumerate}

\item Now there are several very good books on CT for beginners with
  lots of diagrams --- for example \cite{FongSpivak}, \cite{Perrone},
  and \cite{MilewskiCTFPOCaml}. I want to try do draw the ``missing
  diagrams'' for some of their sections, show them to some people, and
  see if they find them useful.

\item I need to learn more Idris and Idris-ct --- and then 1) draw the
  missing diagrams for some of the modules in the Idris-ct sources (as
  a visual guide for the names of the data structures and their
  fields), 2) implement some of my diagrams on Idris-ct; a column with
  Idris-ct code would be a nice addition to, for example, Section
  \ref{basic-example-full}.

\item The paper \cite{PH2} that I uploaded to Arxiv is a kind of
  ``Sheaves for Children'', and some philosopher friends of mine who
  study Alain Badiou --- who uses toposes and sheaves in books like
  \cite{BadiouLoW} and \cite{BadiouMoTT} --- expressed a lot of
  interest in it... the first six sections of \cite{PH2} are
  impeccable (I think!) but the last ones, that are the ones that
  involve categories, were written in a hurry. I need to rewrite them
  using techniques like the ones in Section \ref{teaching-adjunctions}
  to turn them into something like a ``Let's read some sections of
  \cite{Elephant1} and \cite{Riehl} --- an illustrated guide''...
  until I finish that I can't advertise \cite{PH2}, I am too
  embarassed by its last sections.

\end{enumerate}


\newpage

%

\printbibliography
\pu

\end{document}